\documentclass[reqno,12pt]{amsart} 

\usepackage[all]{xy}
\usepackage{supertabular}

\def\theoremnumberstyle{section} 
\def\proofname{Proof}
\usepackage{tom}

\newcommand{\tom}[1]{}

\numberwithin{equation}{section} 
\allowdisplaybreaks[3]

\SelectTips{cm}{12}  
\CompileMatrices     

\begin{document}

   \parindent0cm

     \title[Smoothness]{Smoothness of Equisingular Families of Curves}
     \author{Thomas Keilen}
   \address{Universit\"at Kaiserslautern\\
     Fachbereich Mathematik\\
     Erwin-Schr\"odinger-Stra\ss e\\
     D -- 67663 Kaiserslautern
     }
   \email{keilen@mathematik.uni-kl.de}
   \urladdr{http://www.mathematik.uni-kl.de/\textasciitilde keilen}


   \date{August, 2003.}

     
   \begin{abstract}
     Francesco Severi (cf.~\cite{Sev21}) showed that equisingular
     families of plane nodal curves 
     are T-smooth, i.\ e.\ smooth of the expected dimension, whenever
     they are non-empty. For families with
     more complicated singularities this is no longer true. Given a
     divisor $D$ on a smooth projective surface $\Sigma$ it thus
     makes sense to look for conditions which ensure that the family
     $V_{|D|}^{irr}\big(\ks_1,\ldots,\ks_r\big)$ of irreducible curves
     in the linear system $|D|_l$ with precisely $r$ singular points
     of types $\ks_1,\ldots,\ks_r$ is T-smooth. Considering
     different surfaces including the projective plane, general surfaces in
     $\PC^3$, products of curves and geometrically ruled surfaces, we
     produce a sufficient condition of the type
     \begin{displaymath}
       \sum\limits_{i=1}^r\gamma_\alpha(\ks_i)
       <
       \gamma\cdot (D- K_\Sigma)^2,       
     \end{displaymath}
     where $\gamma_\alpha$ is some invariant of the singularity type
     and $\gamma$ is some constant.
     This generalises the results in \cite{GLS00} for the plane case,
     combining their methods and the method of Bogomolov instability,
     used in \cite{CS97} and \cite{GLS97}. For many singularity types the
     $\gamma_\alpha$-invariant leads to essentially better conditions
     than the invariants used in \cite{GLS97}, and for
     most classes of geometrically ruled
     surfaces our results are the first known  for T-smoothness at
     all. 
   \end{abstract}

   \maketitle



   \section{Introduction}

   The varieties $V_{|D|}(rA_1)$ (respectively the open subvarieties
   $V^{irr}_{|D|}(rA_1)$) 
   of reduced (respectively reduced and irreducible) nodal curves in a fixed
   linear system $|D|_l$ on a smooth projective surface $\Sigma$ are also
   called \emph{Severi varieties}. When $\Sigma=\PC^2$ Severi showed
   that these varieties are smooth of the expected dimension, whenever
   they are non-empty -- that is, nodes always impose independent
   conditions. It seems natural to study this question on other
   surfaces, but it is not surprising that the situation becomes
   harder. 

   Tannenbaum showed in \cite{Tan82} that also on
   K3-surfaces
   $V_{|D|}(rA_1)$ is always smooth, that, however, the dimension is
   larger than the expected one and thus $V_{|D|}(rA_1)$ is not
   T-smooth in this situation. If we  restrict our attention
   to the subvariety $V^{irr}_{|D|}(rA_1)$ of \emph{irreducible} curves
   with $r$ nodes, then we gain T-smoothness again whenever the
   variety is non-empty.
   That is, while on a K3-surface the conditions
   which nodes impose on irreducible curves are always independent,
   they impose dependent conditions on reducible curves. 

   On more complicated surfaces the situation becomes even worse. Chiantini and
   Sernesi study in \cite{CS97} Severi varieties on surfaces
   in $\PC^3$. They show that on a generic quintic $\Sigma$ in $\PC^3$ with
   hyperplane section $H$ the variety
   $V^{irr}_{|dH|}\big(\tfrac{5d(d-2)}{4}\cdot A_1\big)$ has a non-smooth
   reduced component of the expected dimension, if $d$ is even. They
   construct their 
   examples by intersecting a general cone over $\Sigma$ in $\PC^4$
   with a general complete intersection surface of type
   $\big(2,\tfrac{d}{2}\big)$ in $\PC^4$ and projecting the resulting
   curve to $\Sigma$ in $\PC^3$. Moreover, Chiantini and Ciliberto
   give in \cite{CC99} examples showing that the Severi varieties
   $V^{irr}_{|dH|}(rA_1)$ on a
   surface in $\PC^3$ also may have components of dimension larger than the
   expected one. 

   Hence, one can only ask for
   numerical conditions ensuring that $V^{irr}_{|dH|}(rA_1)$ is T-smooth,
   and Chiantini and Sernesi answer this question by showing that on a
   surface of degree $n\geq 5$ the condition 
   \begin{equation}\label{eq:smoothness-intro:1}
     r< \frac{d(d-2n+8)n}{4}
   \end{equation}
   implies that $V^{irr}_{|dH|}(rA_1)$  is
   T-smooth for $d> 2n-8$. Note that the above example shows that this
   bound is even
   sharp. Actually Chiantini and Sernesi prove a somewhat more general result for
   surfaces with ample canonical divisor $K_\Sigma$ and curves which
   are in $|pK_\Sigma|_l$ for some $p\in\Q$. For their proof they
   suppose that for some curve $C\in V^{irr}_{|dH|}(rA_1)$ the
   cohomology group $H^1\big(\Sigma,\kj_{X^*(C)/\Sigma}(D)\big)$
   does not vanish and derive from this the existence of a Bogomolov
   unstable rank-two bundle $E$. This bundle in turn provides them with
   a curve $\Delta$ of small degree realising a large part of the
   zero-dimensional scheme $X^*(C)$, which leads to the desired
   contradiction. 

   This
   is basically the same approach used in \cite{GLS97}. However, they
   allow arbitrary singularities rather than only nodes, and get in
   the case of a surface in $\PC^3$ of degree $n$
   \begin{displaymath}
     \sum_{i=1}^r\big(\tau^*_{ci}(\ks_i)+1\big)^2         
     <  d\cdot\big(d-(n-4)\cdot 
     \max\big\{\tau^*_{ci}(\ks_i)+1\;\big|\;i=1,\ldots,r\big\}\big)\cdot n
   \end{displaymath}
   as main condition for T-smoothness of
   $V_{|dH|}^{irr}(\ks_1,\ldots,\ks_r)$, which for nodal curves
   coincides with 
   \eqref{eq:smoothness-intro:1}. Moreover, for families of plane
   curves of degree $d$ their result gives
   \begin{displaymath}
     \sum_{i=1}^r\big(\tau^*_{ci}(\ks_i)+1\big)^2         
     <  d^2+6d
   \end{displaymath}
   as sufficient condition for T-smoothness, which is weaker than the
   sufficient condition 
   \begin{equation}\label{eq:smoothness-intro:2}
     \sum_{i=1}^r\gamma_1^*(\ks_i)         
     \leq (d+3)^2
   \end{equation}
   derived in \cite{GLS00} using the Castelnuovo function in order to
   provide a curve of small degree which realises a large part of
   $X^*(C)$. The advantage of the $\gamma_1^*$-invariant is that, while
   always bounded from above by $(\tau^*_{ci}+1)^2$, in many
   cases it is substantially smaller -- e.\ g.\ for an
   ordinary $m$-fold point $M_m$, $m\geq 3$, we have $\gamma_1^{es}(M_m)=2m^2$,
   while 
   \begin{displaymath}
     \big(\tau^{es}_{ci}(M_m)+1\big)^2\geq\frac{(m^2+2m+4)^2}{16}.
   \end{displaymath}

   In this paper we combine the methods of \cite{GLS00} and the method of
   Bogomolov instability to reproduce the result
   \eqref{eq:smoothness-intro:2} in the plane case, and to derive
   a similar sufficient condition,
   \begin{displaymath}
     \sum\limits_{i=1}^r\gamma_\alpha(\ks_i)
     <
     \gamma\cdot (D- K_\Sigma)^2,            
   \end{displaymath}
   for T-smoothness on other surfaces --
   involving a generalisation $\gamma_\alpha^*$ of the
   $\gamma_1^*$-invariant which is always bounded from above by the
   latter one. 

   Note that a series 
   of irreducible plane curves of degree $d$ with $r$ singularities of
   type $A_k$, $k$ arbitrarily large, satisfying 
   \begin{displaymath}
     r\cdot k^2=\sum_{i=1}^r\tau^*(A_k)^2=9d^2+\,\text{terms of lower order}
   \end{displaymath}
   constructed by Shustin (cf.~\cite{Shu97}) shows that
   asymptotically\index{asymptotical behaviour} we cannot expect
   to do essentially better in general. For a survey on other known results on
   $\Sigma=\PC^2$ we refer to \cite{GLS00}, and for results on Severi
   varieties on other surfaces see \cite{Tan80,GK89,GLS98b,FM01,Fla02}. 


   In this section we introduce the basic concepts and notations used
   throughout the paper, and we state several important known
   facts. Section \ref{sec:smooth} contains the main results and
   Section \ref{sec:proofs} their
   proofs.


   \subsection{General Assumptions and Notations}\label{subsec:notations}
     Throughout this article $\Sigma$ will denote a smooth projective surface
     over $\C$. 
        
     We will denote by $\Div(\Sigma)$ the 
     group of divisors on $\Sigma$ and by $K_\Sigma$ its canonical
     divisor. 
     If $D$ is any
     divisor on $\Sigma$, $\ko_\Sigma(D)$ shall be the corresponding invertible
     sheaf and we will sometimes write $H^\nu(X,D)$
     instead of $H^\nu\big(X,\ko_X(D)\big)$.  
     A \emph{curve} $C\subset\Sigma$ will be an effective (non-zero) divisor, that
     is a one-dimensional locally principal scheme, not necessarily
     reduced; however, an \emph{irreducible curve} shall be
     reduced by definition.
     $|D|_l$ denotes the
     system of curves linearly equivalent to $D$. 
     We will use the notation
     $\Pic(\Sigma)$ for the \emph{Picard group} of $\Sigma$, that is
     $\Div(\Sigma)$ modulo linear equivalence (denoted by $\sim_l$), and
     $\NS(\Sigma)$ for the 
     \emph{N\'eron--Severi group}, that is $\Div(\Sigma)$ modulo algebraic
     equivalence (denoted by $\sim_a$).
     Given a reduced curve $C\subset\Sigma$ we will write $g(C)$ for
     its \emph{geometric genus}.
   
     Given any closed subscheme
     $X$ of a scheme $Y$, we denote by 
     $\kj_X=\kj_{X/Y}$ the  
     \emph{ideal sheaf} of $X$ in $\ko_Y$. If $X$ is zero-dimensional we denote
     by 
     $\deg(X)=\sum_{z\in Y}\dim_\C(\ko_{Y,z}/\kj_{X/Y,z})$ 
     its \emph{degree}. 
     If $X\subset \Sigma$ is a zero-dimensional
     scheme on $\Sigma$ and $D\in\Div(\Sigma)$, we denote by
     $\big|\kj_{X/\Sigma}(D)\big|_l$ the linear system of curves $C$ in
     $|D|_l$ with $X\subset C$.
     
     Given two curves $C$ and $D$ in $\Sigma$ and a point $z\in
     \Sigma$, and let $f,g\in\ko_{\Sigma,z}$ be local equations at $z$
     of $C$ and $D$ respectively, then we will denote by
     $i(C,D;z)=i(f,g)=\dim_\C(\ko_{\Sigma,z}/\langle f,g\rangle)$ the
     intersection multiplicity of $C$ and $D$ at $z$.

   \subsection{Singularity Types}\label{subsec:types}
     The germ $(C,z)\subset(\Sigma,z)$ of a reduced curve $C\subset\Sigma$ at
     a point $z\in\Sigma$ is called a \emph{plane curve singularity},
     and two plane curve singularities $(C,z)$ and $\big(C',z'\big)$ are said to be
     \emph{topologically} (respectively \emph{analytically equivalent})
     if there is homeomorphism (respectively an analytical isomorphism) 
     $\Phi:(\Sigma,z)\rightarrow(\Sigma,z')$ such that
     $\Phi(C)=C'$. We call an equivalence class with respect to these
     equivalence relations a \emph{topological} (respectively
     \emph{analytical}) \emph{singularity type}. 

     When dealing with numerical conditions for T-smoothness some
     topological (respectively analytical) invariants of the
     singularities play an important role. We gather some results on
     them here for the convenience of the reader.

     Let $(C,z)$ be the germ at $z$ of a reduced curve $C\subset\Sigma$
     and let $f\in R=\ko_{\Sigma,z}$ be a representative of $(C,z)$ in local
     coordinates $x$ and $y$. For the analytical type of the singularity the
     \emph{Tjurina ideal}
     \begin{displaymath}
       I^{ea}(f)=\left\langle\frac{\partial f}{\partial
           x},\frac{\partial f}{\partial y},f \right\rangle
     \end{displaymath}
     plays a very important role, as does the \emph{equisingularity
       ideal} 
     \begin{displaymath}
       I^{es}(f)=\big\{g\in R\;\big|\;f+\varepsilon g \mbox{ is equisingular
         over } \C[\varepsilon]/(\varepsilon^2)\big\}\supseteq I^{ea}(f)
     \end{displaymath}
     for the topological type. They give rise to the following
     invariants of the topological (respectively analytical)
     singularity type $\ks$ of $(C,z)$.

     \begin{enumerate}
     \item  Analytical Invariants:
       \begin{enumerate}
       \item  
         $\tau(\ks)=\dim_\C\big(R/I^{ea}(f)\big)$ is
         the \emph{Tjurina number}, i.\ e.\ the dimension of the base
         space of the semiuniversal deformation of $(C,z)$.
       \item
         $\tau_{ci}(\ks)=\max\big\{\dim_\C(R/I)\;\big|\; 
         I^{ea}(f)\subseteq
         I\mbox{ a complete intersection}\big\}$.
       \item
         $\gamma_\alpha^{ea}(\ks)=\max\big\{\gamma_\alpha(f;I)\;\big|\;
         I^{ea}(f)\subseteq I 
         \mbox{ a complete intersection}\big\}$.
       \end{enumerate}
     \item Topological Invariants:
       \begin{enumerate}
       \item
         $\tau^{es}(\ks)=\dim_\C\big(R/I^{es}(f)\big)$ 
         is the codimension of the $\mu$-constant
         stratum in the semiuniversal deformation of
         $(C,z)$.
       \item
         $\tau^{es}_{ci}(\ks)=\max\big\{\dim_\C(R/I)\;\big|\; 
         I^{es}(C,z)\subseteq I\mbox{ a complete intersection}\big\}$. 
       \item
         $\gamma_\alpha^{es}(\ks)=\max\big\{\gamma_\alpha(f;I)\;\big|\;
         I^{es}(C,z)\subseteq I  
         \mbox{ a complete intersection}\big\}$.
       \end{enumerate}
     \end{enumerate}

     Here, for an ideal $I$ containing $I^{ea}(f)$ and a rational
     number $0\leq \alpha\leq 1$ we define
     \begin{displaymath}
       \gamma_\alpha(f;I)=\max\left\{(1+\alpha)^2\cdot\dim_\C(R/I),\;\lambda_\alpha(f;I,g)\;
         \big|\; g\in I, i(f,g)\leq 2\cdot\dim_\C(R/I)\right\},
     \end{displaymath}
     where for $g\in I$
     \begin{displaymath}
       \lambda_\alpha(f;I,g)=\frac{\big(\alpha\cdot
         i(f,g)-(1-\alpha)\cdot \dim_\C(R/I)\big)^2}{i(f,g)-\dim_\C(R/I)}.
     \end{displaymath}
     Note that by Lemma \ref{lem:schnittzahl} $i(f,g)>\dim_\C(R/I)$
     for all $g\in I$ and $\gamma_\alpha(f,g)$ is thus a well-defined positive rational
     number.

     \bigskip
     \begin{center}
       \framebox[11cm]{
         \begin{minipage}{10cm}
         \medskip
           Throughout this article we will frequently treat
           topological and analytical singularities at the same time.
           Whenever we do so, we will write $\tau^*(\ks)$ for $\tau^{es}(\ks)$
           respectively for $\tau(\ks)$, and analogously we use the
           notation $\tau^*_{ci}(\ks)$ and $\gamma^*_\alpha(\ks)$.
         \medskip
         \end{minipage}
         }
     \end{center}
     \bigskip

     One easily sees the following relations:
     \begin{equation}\label{eq:smoothness-intro:3}
       (1+\alpha)^2\cdot\tau^*_{ci}(\ks)\leq \gamma^*_\alpha(\ks)\leq
       \big(\tau^*_{ci}(\ks)+\alpha\big)^2
       \leq \big(\tau^*(\ks)+\alpha\big)^2.
     \end{equation}

     In \cite{LK03} the $\gamma_\alpha^*$-invariant is calculated for
     the simple singularities,
     \begin{displaymath}
       \begin{array}{|cc|c|}
         \hline
         \multicolumn{2}{|c|}{\ks}  & \gamma_\alpha^{ea}(\ks)= \gamma_\alpha^{es}(\ks)
         \\\hline\hline
         A_k, & k\geq 1 & (k+\alpha)^2 \\\hline
         D_k, & 4\leq k\leq 4+\sqrt{2}\cdot(2+\alpha) 
         & \frac{(k+2\alpha)^2}{2} \\\hline
         D_k, & k\geq 4+\sqrt{2}\cdot(2+\alpha) & (k-2+\alpha)^2\\\hline
         E_k, & k=6,7,8 & \frac{(k+2\alpha)^2}{2}  \\\hline
       \end{array}
     \end{displaymath}
     and for the
     topological singularity type $M_m$ of an ordinary $m$-fold
     point
     \begin{displaymath}
       \gamma_\alpha^{es}(M_m)=2\cdot (m-1+\alpha)^2.
     \end{displaymath}
     Moreover, upper and lower bounds for the
     $\gamma^{es}_0$-invariant and for the $\gamma_1^{es}$-invariant
     of a topological singularity type given by a convenient
     semi-quasihomogeneous power series can be found there. They also
     show that 
     \begin{displaymath}
       \tau_{ci}^{es}(M_m)=
       \left\{
         \begin{array}{cl}
           \frac{(m+1)^2}{4}, & \mbox{ if } m\geq 3 \mbox{ odd},\\
           \frac{m^2+2m}{4}, & \mbox{ if } m\geq 4 \mbox{ even},\\
           1,& \mbox{ if } m=2.
         \end{array}
       \right.
     \end{displaymath}

     These results show in particular that the upper bound for
     $\gamma^*_\alpha(\ks)$ in \eqref{eq:smoothness-intro:3} may be
     attained, while it may as well be far from the actual value.

     \begin{lemma}\label{lem:schnittzahl}
       Let $(C,z)$ be a reduced plane curve singularity given by
       $f\in\ko_{\Sigma,z}$ and let
       $I\subseteq\m_{\Sigma,z}\subset\ko_{\Sigma,z}$ be an ideal
       containing the Tjurina ideal $I^{ea}(C,z)$. Then for any $g\in I$
       we have
       \begin{displaymath}
         \dim_\C(\ko_{\Sigma,z}/I)<\dim_\C\big(\ko_{\Sigma,z}/(f,g)\big)=i(f,g).
       \end{displaymath}
     \end{lemma}
     \begin{proof}
       Cf.~\cite{Shu97} Lemma 4.1.
     \end{proof}

   \subsection{Singularity Schemes}\label{subsec:schemes}          

     For a reduced curve $C\subset\Sigma$ we recall the definition of
     the zero-dimensional schemes $X^{es}(C)$ and
     $X^{ea}(C)$ from \cite{GLS00}. They are defined
     by the ideal sheaves $\kj_{X^{es}(C)/\Sigma}$ and
     $\kj_{X^{ea}(C)/\Sigma}$
     respectively, given by the stalks
     $\kj_{X^{es}(C)/\Sigma,z}=I^{es}(f)$ and
     $\kj_{X^{ea}(C)/\Sigma,z}=I^{ea}(f)$ respectively,
     where $f\in\ko_{\Sigma,z}$ is a local equation of $C$ at $z$.
     We call $X^{es}(C)$ the \emph{equisingularity scheme} of $C$ 
     and
     $X^{ea}(C)$ the \emph{equianalytical singularity scheme} of
     $C$. 

     \bigskip
     \begin{center}
       \framebox[11cm]{
         \begin{minipage}{10cm}
         \medskip
           Throughout this article we will frequently treat
           topological and analytical singularities at the same time.
           Whenever we do so, we will write $X^*(C)$ for $X^{es}(C)$
           respectively for $X^{ea}(C)$.
         \medskip
         \end{minipage}
         }
     \end{center}
     \bigskip

   \subsection{Equisingular Families}\label{subsec:families}
     Given a divisor $D\in\Div(\Sigma)$ and topological or analytical singularity types
     $\ks_1,\ldots, \ks_r$, we denote
     by $V=V_{|D|}(\ks_1,\ldots,\ks_r)$ the locally closed subspace of
     $|D|_l$ of reduced curves in the linear system $|D|_l$ having
     precisely $r$ singular points of types $\ks_1,\ldots,\ks_r$.      
     By
     $V^{irr}=V_{|D|}^{irr}(\ks_1,\ldots,\ks_r)$ we denote the
     open  subset of $V$ of irreducible curves.  
     If a type $\ks$ occurs $k>1$ times, we rather write $k\ks$
     than $\ks,\stackrel{k}{\ldots},\ks$.
     We call these families of curves \emph{equisingular families of
       curves}.

     We say that $V$ is \emph{T-smooth} at $C\in V$ if the germ
     $(V,C)$ is smooth of the (expected) dimension
     $\dim|D|_l-\deg\big(X^*(C)\big)$.
     By \cite{Los98}
     Proposition 2.1 (see also \cite{GK89}, \cite{GL96}, \cite{GLS00})
     T-smoothness of $V$ at $C$ follows from the vanishing of
     $H^1\big(\Sigma,\kj_{X^*(C)/\Sigma}(C)\big)$, since
     the tangent space of $V$ at $C$ may be
     identified with
     $H^0\big(\Sigma,\kj_{X^*(C)/\Sigma}(C)\big)/H^0(\Sigma,\ko_\Sigma)$.


   \section{The Main Results}\label{sec:smooth}

   In this section we give sufficient conditions for the
   T-smoothness of equisingular families of curves on certain
   surfaces with Picard number one, including the projective
   plane, general surfaces in $\PC^3$ and general K3-surfaces --, on
   general products of curves, and on  geometrically ruled
   surfaces.

   \subsection{Surfaces with Picard Number One}
   
   \begin{theorem}\label{thm:smooth-add-p3}
     Let $\Sigma$ be a surface such that
     $\NS(\Sigma)=L\cdot\Z$ with $L$ ample, 
     let $D=d\cdot L\in\Div(\Sigma)$, let $\ks_1,\ldots,\ks_r$ be 
     topological or analytical singularity types, and let
     $K_\Sigma=\kappa\cdot L$.
     Suppose that $d\geq\max\{\kappa+1,-\kappa\}$ and
     \begin{equation}
       \label{eq:smooth-add-p3:1}
       \sum\limits_{i=1}^r\gamma_\alpha^*(\ks_i)<\alpha\cdot
       (D-K_\Sigma)^2=\alpha\cdot (d-\kappa)^2\cdot L^2\;\;\;
       \mbox{ with } \alpha= \tfrac{1}{\max\{1,1+\kappa\}}.
     \end{equation}     
     Then either
     $V_{|D|}^{irr}(\ks_1,\ldots,\ks_r)$ is empty
     or it is T-smooth.
     \hfill $\Box$
   \end{theorem}



   We now apply the result in several special cases.

   \begin{corollary}\label{cor:smooth-add-p2}
     Let $d\geq 3$, $L\subset\PC^2$ be a line, and 
     $\ks_1,\ldots,\ks_r$ be 
     topological or analytical singularity types.
     Suppose that 
     \begin{equation}\label{eq:smooth-add-p2:1}
       \sum\limits_{i=1}^r\gamma_1^*(\ks_i) <
       (d+3)^2.
     \end{equation}     
     Then either $V_{|dH|}^{irr}(\ks_1,\ldots,\ks_r)$ is empty or
     T-smooth.
     \hfill $\Box$
   \end{corollary}

   As soon as for one of the singularities we have
   $\gamma^*_1(\ks_i)>4\cdot\tau^*_{ci}(\ks_i)$, e.\ g.\ simple singularities
   or ordinary multiple points which are not simple double points,
   then the strict inequality in
   (\ref{eq:smooth-add-p2:1}) can be replaced by
   ``$\leq$'', which then is the same sufficient condition as 
   in \cite{GLS00} Theorem 1 (see also \eqref{eq:smoothness-intro:2}).

   In particular, 
   $V_{|dH|}^{irr}(kA_1,mA_2,M_{m_1},\ldots,M_{m_r})$, $m_i\geq 3$,  is therefore T-smooth as
   soon as
   \begin{displaymath}
     4k+9m+\sum_{i=1}^r 2\cdot m_i^2\leq (d+3)^2.
   \end{displaymath}

   For further results in
   the plane case see
   \cite{Wah74a,GK89,Lue87,Lue87a,Shu87,Vas90,Shu91a,Shu94,GL96,Shu96a,Shu97,GLS98b,Los98,GLS00}.

   A smooth complete intersection surface with Picard
   number one satisfies the assumptions of Theorem
   \ref{thm:smooth-add-p3}. Thus by the Theorem of Noether the result
   applies in particular to general surfaces in $\PC^3$.
   Moreover, if in Theorem \ref{thm:smooth-add-p3} we have $\kappa>0$, i.\ e.\
   $\alpha<1$, then
   the strict inequality in Condition \eqref{eq:smooth-add-p3:1} may
   be replaced by ``$\leq$'', since in \eqref{eq:smooth-add-p3:2} the second
   inequality is strict, as is the second inequality in \eqref{eq:smooth-add-p3:3}.

   \begin{corollary}
     Let $\Sigma\subset\PC^3$ be a smooth hypersurface of degree
     $n\geq 4$,
     let $H\subset\Sigma$ be a hyperplane
     section, and suppose that the Picard number  of
     $\Sigma$ is one.
     Let $d\geq n-3$ and let $\ks_1,\ldots,\ks_r$ be 
     topological or analytical singularity types.
     Suppose that 
     \begin{displaymath}
       \sum\limits_{i=1}^r\gamma_{\frac{1}{n-3}}^*(\ks_i)\leq\frac{n}{n-3}\cdot (d-n+4)^2.
     \end{displaymath}     
     Then either
     $V_{|D|}^{irr}(\ks_1,\ldots,\ks_r)$ is empty
     or it is T-smooth.
     \hfill $\Box$
   \end{corollary}

   In particular, 
   $V_{|dH|}^{irr}(M_{m_1},\ldots,M_{m_r})$, $m_i\geq 3$,  is therefore T-smooth as
   soon as
   \begin{displaymath}
     \sum_{i=1}^r 2\cdot \left(m_i-\frac{n-2}{n-3}\right)^2< \frac{n}{n-3}\cdot (d-n+4)^2,
   \end{displaymath}
   which is better than the conditions derived from \cite{GLS97}.
   The condition 
   \begin{displaymath}
     r\leq \frac{n\cdot (n-3)}{(n-2)^2}\cdot (d-n+4)^2,
   \end{displaymath}
   which gives the T-smoothness of $V_{|dH|}(rA_1)$ is weaker than the
   condition provided in \cite{CS97}, but for
   $n=5$ it reads
   $r\leq \frac{10}{9}\cdot (d-1)^2$
   and comes still  close to the sharp bound
   $\frac{5}{4}\cdot (d-1)^2$ provided there for odd $d$.

   A general K3-surface has also Picard number one..

   \begin{corollary}
     Let $\Sigma$ be a smooth K3-surface with $\NS(\Sigma)=L\cdot\Z$,
     $L$ ample, and set $n=L^2$.
     Let $d\geq 1$, and let  $\ks_1,\ldots,\ks_r$ be 
     topological or analytical singularity types.
     Suppose that
     \begin{displaymath}
       \sum\limits_{i=1}^r\gamma_1^*(\ks_i)         
       <  d^2n.
     \end{displaymath}
     Then either
     $V_{|dL|}^{irr}(\ks_1,\ldots,\ks_r)$ is empty
     or it is T-smooth.
     \hfill $\Box$
   \end{corollary}

   The best previously known condition for T-smoothness on K3-surfaces
   \begin{displaymath}
     \sum\limits_{i=1}^r\big(\tau^*_{ci}(\ks_i)+1\big)^2         
       <  d^2n
   \end{displaymath}
   is thus completely replaced. 


   \subsection{Products of Curves}\label{subsec:smooth-add:products-of-surfaces}

   If $\Sigma=C_1\times C_2$ is the product of two smooth projective
   curves, then for a general choice of $C_1$ and $C_2$ the
   N\'eron--Severi group will be generated by two fibres of the
   canonical projections, by abuse of notation also denoted by $C_1$
   and $C_2$. If both curves are elliptic, then ``general'' just means
   that the two curves are non-isogenous.

   \begin{theorem}\label{thm:smooth-add-products-of-curves}
     Let $C_1$ and $C_2$ be two smooth projective curves of genera $g_1$
     and $g_2$ with $g_1\geq g_2$, such that for $\Sigma=C_1\times C_2$
     the N\'eron--Severi group is $\NS(\Sigma)=C_1\Z\oplus C_2\Z$.

     Let $D\in\Div(\Sigma)$ such that $D\sim_a aC_1+bC_2$ with
     $a\geq\max\big\{2-2g_2,2g_2-1\big\}$ 
     and
     $b\geq\max\big\{2-2g_1,2g_1-1\big\}$, let
     $\ks_1,\ldots,\ks_r$ be 
     topological or analytical singularity types.
     Suppose that 
     \begin{equation}\label{eq:smooth-add-products-of-curves:1}
       \sum\limits_{i=1}^r\gamma_0^*(\ks_i) 
       \;<\;\gamma\cdot(D-K_\Sigma)^2, 
     \end{equation}
     where the constant $\gamma$ may be read off  the following
     table with
     $A=\frac{a-2g2+2}{b-2g1+2}$
     \begin{displaymath}
       \begin{array}{|r|r|c|}
         \hline
         g_1 & g_2 & \gamma \\\hline
         0,1 & 0,1 & \frac{1}{4} \\\hline
         \geq 2 & 0,1 & \min\left\{\frac{1}{4g_1},
           \frac{1}{4\cdot(g_1-1)\cdot A}\right\}\\\hline
         \geq 2 & \geq 2 & \min\left\{\frac{1}{4g_1+4g_2-4},
           \frac{A}{4\cdot(g_2-1)},\frac{1}{4\cdot(g_1-1)\cdot A}\right\}\\\hline
       \end{array}
     \end{displaymath} 
     Then either
     $V_{|D|}^{irr}(\ks_1,\ldots,\ks_r)$ is empty
     or it is T-smooth.
     \hfill$\Box$
   \end{theorem}

   In particular, on a product of non-isogenous elliptic curves for
   nodal curves we reproduce the previous sufficient condition 
   \begin{displaymath}
     r<\frac{ab}{2},
   \end{displaymath}
   for
   T-smoothness of $V_{|aC_1+bC_2|}^{irr}(rA_1)$ from \cite{GLS97}, 
   while the
   previous general condition
   \begin{displaymath}
     \frac{\big(m_i^2+2m_i+5\big)^2}{32}<ab
   \end{displaymath}
   for T-smoothness of 
   $V_{|aC_1+bC_2|}^{irr}(M_{m_1},\ldots,M_{m_r})$, $m_i\geq 3$,   
   has been replaced by 
   \begin{displaymath}
     \sum_{i=1}^r 4\cdot (m_i-1)^2< ab,
   \end{displaymath}
   which is better from $m_i=7$ on.
   
   Note that the constant $\gamma$ in Theorem
   \ref{thm:smooth-add-products-of-curves} depends on the ratio of $a$
   and $b$ unless both $g_1$ and $g_2$ are at most one. 
   This means
   that in general an asymptotical behaviour can
   only be examined if the ratio remains unchanged.


   \subsection{Geometrically Ruled Surfaces}\label{subsec:smooth-add:ruled-surfaces}

   Let $\pi:\Sigma=\PC(\ke)\rightarrow C$ be a geometrically ruled
   surface with normalised bundle $\ke$ (in the 
   sense of \cite{Har77} V.2.8.1). The N\'eron--Severi group of
   $\Sigma$ is $\NS(\Sigma) = C_0\Z\oplus F\Z$ 
   with intersection matrix $\left(\begin{smallmatrix}-e & 1 \\ 1 & 0\end{smallmatrix}\right)$
   where $F\cong\PC^1$ is a fibre of $\pi$, $C_0$ a section of $\pi$
   with $\ko_\Sigma(C_0)\cong\ko_{{\mathbbm P}(\ke)}(1)$, $g=g(C)$ the genus of
   $C$,  $\mathfrak{e}=\Lambda^2\ke$ and
   $e=-\deg(\mathfrak{e})\geq -g$. 
   For the canonical divisor we have $K_\Sigma \sim_a -2C_0+ (2g-2-e)\cdot F$.

   \begin{theorem}\label{thm:smooth-add-ruled-surfaces}
     Let $\pi:\Sigma\rightarrow C$ be a geometrically ruled surface
     with  $g=g(C)$.
     Let $D\in\Div(\Sigma)$ such that $D\sim_a aC_0+bF$ with
     $b>\max\{2g-2,2-2g\}+\frac{ae}{2}$ and $a> 2$, and
     let $\ks_1,\ldots,\ks_r$ be 
     topological or analytical singularity types.
     Suppose that 
     \begin{equation}\label{eq:smooth-add-ruled-surfaces:1}
       \sum\limits_{i=1}^r\gamma_0^*(\ks_i)
       \;<\;\gamma\cdot(D-K_\Sigma)^2,
     \end{equation} 
     where with
     $A=\frac{a+2}{b+2-2g-\frac{ae}{2}}$ the constant $\gamma$ satisfies
     \begin{displaymath}
       \gamma=
       \left\{
         \begin{array}{ll}
           \frac{1}{4}, & \mbox{ if } g\in\{0,1\},\\
           \min\left\{\frac{1}{4g},\frac{1}{4\cdot(g-1)\cdot
               A}\right\},& \mbox{ if } g\geq 2.
         \end{array}
       \right.
     \end{displaymath} 
     Then either
     $V_{|D|}^{irr}(\ks_1,\ldots,\ks_r)$ is empty
     or it is T-smooth.
     \hfill $\Box$
   \end{theorem}

   The results of \cite{GLS97} only applied to eight Hirzebruch
   surfaces and a few classes of fibrations over elliptic curves,
   while our results apply to all geometrically ruled
   surfaces. Moreover, the results are in general better, e.\ g.\ for
   the Hirzebruch surface $\PC^1\times\PC^1$ already the previous
   sufficient condition for T-smoothness of families of curves with
   $r$ cusps and $b=3a$ the condition 
   \begin{displaymath}
     9r<2a^2+8a
   \end{displaymath}
   has been replaced by the slightly better condition
   \begin{displaymath}
     8r<3a^2+8a+4.
   \end{displaymath}
   For ordinary multiple points the difference will become more
   significant. Even for families of nodal curves the new conditions
   would always be slightly better, 
   but for those families T-smoothness is guaranteed anyway by 
   \cite{Tan80}. 


   Note that, as for products of curves, the constant $\gamma$ in Theorem
   \ref{thm:smooth-add-ruled-surfaces} depends on the ratio of $a$
   and $b$ unless $g$ is at most one.


   \section{The Proofs}\label{sec:proofs}

   The following Lemma is the technical key to the above
   results. Using the method of Bogomolov unstable vector bundles, it
   gives us a ``small'' curve which passes through a ``large'' part of
   $X^*(C)$, provided that
   $h^1\big(\Sigma,\kj_{X^*(C)/\Sigma}(D)\big)\not=0$. We will then
   show that its existence contradicts \eqref{eq:smooth-add-p3:1},
   \eqref{eq:smooth-add-products-of-curves:1}, or 
   \eqref{eq:smooth-add-ruled-surfaces:1} respectively.  

   \begin{lemma}\label{lem:smooth-add-A}
     Let $\Sigma$ a smooth projective surface,   
     and let $D\in\Div(\Sigma)$
       and $X\subset \Sigma$  be a
     zero-dimensional scheme satisfying
     \begin{enumerate}
     \item[(0)] $D-K_\Sigma$ is big and nef,  and $D+K_\Sigma$ is nef,
     \item[(1)] $\exists\;C\in|D|_l\;\mbox{ irreducible}:\;X\subseteq X^*(C)$,
     \item[(2)] $h^1\big(\Sigma,\kj_{X/\Sigma}(D)\big)>0$, and 
     \item[(3)] $4\cdot \deg(X_0)< (D-K_\Sigma)^2$ for all
       local complete intersection schemes
       $X_0\subseteq X$.
     \end{enumerate}
     Then there exists a curve $\Delta\subset\Sigma$ and a
     zero-dimensional local complete intersection scheme $X_0\subseteq
     X\cap\Delta$  such that with the notation
     $\supp(X_0)=\{z_1,\ldots,z_s\}$, $X_i=X_{0,z_i}$ 
     and\footnote{Since $X_0\subseteq X^*(C)\subseteq X^{ea}(C)$, Lemma
       \ref{lem:schnittzahl} applies
       to the local ideals of $X_0$, that is for the
       points $z\in \supp(X_0)$ we have
       $i(C,\Delta;z)\geq\deg(X_0,z)+1$. 
       }
     $\varepsilon_i=\min\{\deg(X_i),i(C,\Delta;z_i)-\deg(X_i)\}\geq 1$ we have 
     \begin{enumerate}
     \item $D.\Delta\geq\deg(X_0)+\sum_{i=1}^s\varepsilon_i$,
     \item
       $\deg\big(X_0\big)
       \geq\big(D-K_\Sigma-\Delta\big).\Delta$, 
     \item $\big(D-K_\Sigma-2\cdot\Delta\big)^2>0$, and
     \item  $\big(D-K_\Sigma-2\cdot\Delta\big).H>0$\;\;
       for all $H\in\Div(\Sigma)$ ample.
     \end{enumerate}
     Moreover, it follows
     \begin{equation}
       \label{eq:smooth-add-A:1}
       0\leq \tfrac{1}{4}\cdot(D-K_\Sigma)^2-\deg\big(X_0\big)
       \leq \left(\tfrac{1}{2}\cdot(D-K_\Sigma)-\Delta\right)^2.
     \end{equation}
   \end{lemma}

   \begin{proof}
     Choose $X_0\subseteq X$ minimal such that still
     $h^1\big(\Sigma,\kj_{X_0/\Sigma}(D)\big)>0$.  
     By Assumption (0) the divisor $D-K_\Sigma$ is big and
     nef, and thus $h^1\big(\Sigma,\ko_\Sigma(D)\big)=0$ 
     by the Kawamata--Viehweg Vanishing Theorem. Hence $X_0$ cannot be
     empty. 

     Due to the Grothendieck-Serre duality \tom{(cf.~\cite{Har77}
       III.7.6) }
     we have
     \begin{displaymath}
       0\not=H^1\big(\Sigma,\kj_{X_0/\Sigma}(D)\big)
       \cong
       \Ext^1\big(\kj_{X_0/\Sigma}(D-K_\Sigma),\ko_\Sigma\big).
     \end{displaymath}  
     That is, there is an extension  \tom{(cf.~\cite{Har77}
       Ex.~III.6.1) }
     \begin{equation}\label{eq:smooth-add-A:2}
       0\rightarrow \ko_\Sigma\rightarrow
       E\rightarrow\kj_{X_0/\Sigma}(D-K_\Sigma) 
       \rightarrow 0.
     \end{equation}
     The minimality of $X_0$ implies that $E$ is locally free and
     $X_0$ is a local complete intersection scheme
     (cf.~\cite{Laz97} Proposition 3.9). Moreover, we have
     \tom{(cf.~\cite{Laz97} Exercise 4.3)} 
     \begin{equation}\label{eq:smooth-add-A:3}
       c_1(E)=D-K_\Sigma\mbox{\;\;\;and\;\;\;}
       c_2(E)=\deg(X_0).
     \end{equation}

     By Assumption (3) and  \eqref{eq:smooth-add-A:3} we have 
     \begin{displaymath}
       c_1(E)^2-4\cdot c_2(E)=(D-K_\Sigma)^2-4\cdot \deg(X_0) >0,
     \end{displaymath}
     and thus $E$ is Bogomolov\index{Bogomolov} unstable (cf.~\cite{Laz97} Theorem
     4.2).      
     This, however, implies that there exists a divisor $\Delta_0\in\Div(\Sigma)$ and
     a zero-dimensional scheme $Z\subset\Sigma$  such that
     \begin{equation}
       \label{eq:smooth-add-A:4}
       0\rightarrow\ko_\Sigma(\Delta_0)\rightarrow
       E\rightarrow\kj_{Z/\Sigma}(D-K_\Sigma-\Delta_0)
       \rightarrow 0
     \end{equation}
     is exact\tom{ (cf.~
     \cite{Laz97} Theorem 4.2)}, and such that 
     \begin{equation}
       \label{eq:smooth-add-A:5}
       (2\Delta_0-D+K_{\Sigma})^2\geq c_1(E)^2-4\cdot c_2(E)>0
     \end{equation}
     and
     \begin{equation}
       \label{eq:smooth-add-A:6}
       (2\Delta_0-D+K_{\Sigma}).H> 0 \;\;\;\mbox{ for all ample }\;\;\;
       H\in\Div(\Sigma).
     \end{equation}
     
     Tensoring \eqref{eq:smooth-add-A:4} with
     $\ko_\Sigma(-\Delta_0)$ leads to the following exact sequence
     \begin{equation}\label{eq:smooth-add-A:7}
       0\rightarrow\ko_\Sigma\rightarrow
       E(-\Delta_0)
       \rightarrow\kj_{Z/\Sigma}\left(D-K_\Sigma-2\Delta_0\right)
       \rightarrow 0,
     \end{equation}
     and we deduce 
     $h^0\big(\Sigma,E(-\Delta_0)\big)\not=0$. 

     Now tensoring \eqref{eq:smooth-add-A:2} with $\ko_\Sigma(-\Delta_0)$ leads to
     \begin{equation}\label{eq:smooth-add-A:8}
       0\rightarrow \ko_\Sigma(-\Delta_0)\rightarrow
       E(-\Delta_0)
       \rightarrow\kj_{X_0/\Sigma}\left(D-K_\Sigma-\Delta_0\right) 
       \rightarrow 0.
     \end{equation}
     Let $H$ be some ample divisor. By \eqref{eq:smooth-add-A:6} and since
     $D-K_\Sigma$ is nef by (0):
     \begin{displaymath}
       -\Delta_0.H<-\tfrac{1}{2}\cdot(D-K_\Sigma).H\leq 0.
     \end{displaymath}
     Hence
     $-\Delta_0$ cannot be effective, that is
     $H^0\big(\Sigma,\ko_\Sigma(-\Delta_0)\big)=0$. But the long exact
     cohomology sequence of \eqref{eq:smooth-add-A:8} then implies
     \begin{displaymath}
       0\not=H^0\big(\Sigma,E(-\Delta_0)\big)
       \hookrightarrow 
       H^0\left(\Sigma,\kj_{X_0/\Sigma}
         \left(D-K_\Sigma-\Delta_0\right)\right).
     \end{displaymath}
     In particular 
     we may choose
     a curve 
     \begin{displaymath}
       \Delta\in
       \big|\kj_{X_0/\Sigma}(D-K_\Sigma-\Delta_0)\big|_l.
     \end{displaymath}
     Thus (c) and (d) follow from \eqref{eq:smooth-add-A:5} and
     \eqref{eq:smooth-add-A:6}. It remains to show 
     (a) and (b). 

     We note that $C\in|D|_l$ is irreducible and that
     $\Delta$ cannot contain $C$ as an irreducible component:
     otherwise applying \eqref{eq:smooth-add-A:6} with some ample
     divisor $H$ we would get 
     the following contradiction, since $D+K_\Sigma$ is nef by (0), 
     \begin{displaymath}
       0\leq (\Delta-C).H<-\tfrac{1}{2}\cdot (D+K_\Sigma).H\leq 0.
     \end{displaymath}
     Since $X_0\subset C\cap\Delta$ 
     the Theorem of B\'ezout implies (a):
     \begin{displaymath}
       D.\Delta=C.\Delta=\sum_{z\in C\cap\Delta}i(C,\Delta;z)\geq
       \sum_{i=1}^s\big(\deg(X_i)+\varepsilon_i\big)=\deg(X_0)+\sum_{i=1}^s\varepsilon_i.
     \end{displaymath}
     Finally, by \eqref{eq:smooth-add-A:3} and \eqref{eq:smooth-add-A:4} we get (b):
     \begin{displaymath}
       \deg(X_0)=c_2(E)=\Delta_0.(D-K_\Sigma-\Delta_0)+\deg(Z)
       \geq (D-K_\Sigma-\Delta).\Delta.
     \end{displaymath}
     Equation \eqref{eq:smooth-add-A:1} is just a reformulation of
     (b).
   \end{proof}

   Using this result we can now prove the main theorems.

   \begin{proof}[Proof of Theorem \ref{thm:smooth-add-p3}]
     Let $C\in V_{|D|}^{irr}(\ks_1,\ldots,\ks_r)$. It suffices to show
     that the cohomology group
     $h^1\big(\Sigma,\kj_{X^*(C)/\Sigma}(D)\big)$ vanishes. 

     Suppose this is not the case. Since for $X_0\subseteq X^*(C)$ any
     local complete intersection scheme and
     $z\in\supp(X_0)$ we have
     \begin{equation}\label{eq:smooth-add-p3:2}
       4\cdot\deg(X_z)\leq\frac{4}{(1+\alpha)^2}\cdot\gamma_\alpha^*(C,z)
       \leq\frac{1}{\alpha} \cdot\gamma_\alpha^*(C,z)
     \end{equation}
     Lemma \ref{lem:smooth-add-A} applies and there is curve
     $\Delta\in|\delta\cdot L|_l$ and a local complete intersection
     scheme $X_0\subseteq X^*(C)$ satisfying the assumptions (a)-(d)
     there and Equation~\eqref{eq:smooth-add-A:1}. That is, fixing the
     notation $l=\sqrt{L^2}$, $\supp(X_0)=\{z_1,\ldots,z_s\}$, $X_i=X_{0,z_i}$ and
     $\varepsilon_i=\min\{\deg(X_i),i(C,\Delta;z_i)-\deg(X_i)\}\geq 1$, we
     have  
     \begin{enumerate}
     \item $d\cdot\delta\cdot l^2\geq \deg(X_0)+\sum_{i=1}^s\varepsilon_i$,
     \item $\deg(X_0)\geq (d-\kappa-\delta)\cdot \delta\cdot l^2$,
     \end{enumerate}
     and 
     \begin{displaymath}
       \delta\cdot l\leq \tfrac{(d-\kappa)\cdot l}{2}
       -\sqrt{\tfrac{(d-\kappa)^2\cdot l^2}{4}-\deg(X_0)}
       =
       \frac{2\cdot\deg(X_0)}{
         (d-\kappa)\cdot l
         +\sqrt{(d-\kappa)^2\cdot l^2-4\cdot\deg(X_0)}}.
     \end{displaymath}
     But then together with (a) and (b) we deduce
     \begin{equation}\label{eq:smooth-add-p3:3}
       \sum_{i=1}^s\varepsilon_i\leq \delta\cdot(\delta+\kappa)\cdot l^2
       \leq
       \frac{1}{\alpha} \cdot
       \left(
         \frac{2\cdot\deg(X_0)}{
           (d-\kappa)\cdot l
           +\sqrt{(d-\kappa)^2\cdot l^2-4\cdot\deg(X_0)}}           
       \right)^2. 
     \end{equation}
     Applying the Cauchy inequality this leads to
     \begin{displaymath}
       \sum_{i=1}^s\frac{\deg(X_i)^2}{\varepsilon_i}
       \geq
       \frac{\deg(X_0)^2}{\sum_{i=1}^s\varepsilon_i}
       \geq
       \frac{\alpha\cdot(d-\kappa)^2\cdot l^2}{4}\cdot 
       \left(1+\sqrt{1-\tfrac{4\cdot\deg(X_0)}{(d-\kappa)^2\cdot l^2}}\right)^2.           
     \end{displaymath}
     Setting
     \begin{displaymath}
       \beta=\frac{\sum_{i=1}^s\frac{\deg(X_i)^2}{\varepsilon_i}}{\alpha\cdot(d-\kappa)^2\cdot l^2},
       \;\;\;\;
       \gamma=\frac{\sum_{i=1}^s\frac{\deg(X_i)^2}{\varepsilon_i}}{\alpha\cdot\deg(X_0)},
     \end{displaymath}
     we thus have 
     \begin{displaymath}
       \beta\geq\frac{1}{4}\cdot\left(1+\sqrt{1-\tfrac{4\beta}{\gamma}}\right)^2,
     \end{displaymath}
     and hence, $\beta\geq\big(\tfrac{\gamma}{\gamma+1}\big)^2$.
     But then, applying the Cauchy inequality once more, we find
     \begin{displaymath}
       \alpha\cdot(d-\kappa)^2\cdot l^2
       = 
       \frac{\alpha\cdot\gamma}{\beta}\cdot \deg(X_0)
       \leq 
       \alpha\cdot\left(\gamma+2+\frac{1}{\gamma}\right)\cdot \deg(X_0)
     \end{displaymath}
     \begin{displaymath}
       \leq
       \sum_{i=1}^s\left(\frac{\deg(X_i)^2}{\varepsilon_i}+2\alpha\deg(X_i)+\alpha^2\varepsilon_i\right)
       \leq
       \sum_{i=1}^r\gamma_\alpha^*(\ks_i),
     \end{displaymath}
     in contradiction to Equation \eqref{eq:smooth-add-p3:1}.
   \end{proof}

   \begin{proof}[Proof of Theorem \ref{thm:smooth-add-products-of-curves}]
     Let $C\in V_{|D|}^{irr}(\ks_1,\ldots,\ks_r)$. It suffices to show
     that the cohomology group
     $h^1\big(\Sigma,\kj_{X^*(C)/\Sigma}(D)\big)$ vanishes. 

     Suppose this is not the case. Since for $X_0\subseteq X^*(C)$ any
     local complete intersection scheme and
     $z\in\supp(X)$ we have
     \begin{displaymath}
       \deg(X_z)\leq\gamma_0^*(C,z),
     \end{displaymath}
     and since $\gamma\leq \frac{1}{4}$, 
     Lemma \ref{lem:smooth-add-A} applies and there is curve
     $\Delta\sim_a\alpha\cdot C_1+\beta\cdot C_2$ and a local complete intersection
     scheme $X_0\subseteq X^*(C)$ satisfying the assumptions (a)-(d) there
     and Equation~\eqref{eq:smooth-add-A:1}. 
     That is, fixing the
     notation $\supp(X_0)=\{z_1,\ldots,z_s\}$, $X_i=X_{0,z_i}$ and
     $\varepsilon_i=\min\{\deg(X_i),i(C,\Delta;z_i)-\deg(X_i)\}\geq 1$, we
     have  
     \begin{enumerate}
     \item $a\beta+b\alpha\geq \deg(X_0)+\sum_{i=1}^s\varepsilon_i$,
     \item $\deg(X_0)\geq (a-2g_2+2-\alpha)\cdot \beta+
       (b-2g_1+2-\beta)\cdot \alpha$, and 
     \item $0\leq\alpha\leq \frac{a-2g_2+2}{2}$ and $0\leq\beta\leq \frac{b-2g_1+2}{2}$.
     \end{enumerate}
     The last inequalities follow from (d) in Lemma
     \ref{lem:smooth-add-A} replacing the ample divisor $H$ by the nef
     divisors $C_2$ respectively $C_1$. 

     From (b) and (c) we deduce
     \begin{displaymath}
       \deg(X_0)\geq \frac{a-2g_2+2}{2}\cdot \beta+
       \frac{b-2g_1+2}{2}\cdot \alpha,
     \end{displaymath}
     and thus
     \begin{equation}
       \label{eq:smooth-add-products-of-curves:2}
       \deg(X_0)^2\geq 4\cdot
       \frac{a-2g_2+2}{2}\cdot\frac{b-2g_1+2}{2}\cdot \alpha\cdot\beta
       = \frac{(D-K_\Sigma)^2}{2}\cdot \alpha\cdot \beta.
     \end{equation}
     Considering now (a) and (b) we get
     \begin{displaymath}
       0<\sum_{i=1}^s\varepsilon_i\leq \Delta.(\Delta+K_\Sigma)=
       2\alpha\beta+(2g_1-2)\cdot \alpha + (2g_2-2)\cdot \beta
       \leq \frac{\alpha\beta}{2\gamma},
     \end{displaymath}
     where the last inequality holds only if
     $\alpha\not=0\not=\beta$. In particular, we see
     $\alpha\not=0$ if $g_2\leq 1$ and $\beta\not=0$ if $g_1\leq 1$.
     But this together with \eqref{eq:smooth-add-products-of-curves:2}
     gives
     \begin{displaymath}
       \sum_{i=1}^s\varepsilon_i\leq\frac{\deg(X_0)^2}{\gamma\cdot(D-K_\Sigma)^2}.
     \end{displaymath}

     If $\alpha=0$, then from (a) and (b) we
     deduce again
     \begin{displaymath}
       0<\sum_{i=1}^s\varepsilon_i\leq (2g_2-2)\cdot \beta\leq
       \frac{4\cdot(g_1-1)}{A}\cdot \frac{\deg(X_0)^2}{(D-K_\Sigma)^2}
       \leq \frac{\deg(X_0)^2}{\gamma\cdot(D-K_\Sigma)^2},
     \end{displaymath}
     and similarly, if $\beta=0$,
     \begin{displaymath}
       0<\sum_{i=1}^s\varepsilon_i\leq (2g_1-2)\cdot \alpha\leq
       4\cdot(g_1-1)\cdot A\cdot \frac{\deg(X_0)^2}{(D-K_\Sigma)^2}
       \leq \frac{\deg(X_0)^2}{\gamma\cdot(D-K_\Sigma)^2}.
     \end{displaymath}

     Applying the Cauchy inequality, we finally get
     \begin{displaymath}
       \gamma\cdot(D-K_\Sigma)^2\leq
       \frac{\deg(X_0)^2}{\sum_{i=1}^s\varepsilon_i}
       \leq \sum_{i=1}^s\frac{\deg(X_i)^2}{\varepsilon_i}
       \leq \sum_{i=1}^r\gamma_0^*(\ks_i),
     \end{displaymath}
     in contradiction to Assumption \eqref{eq:smooth-add-products-of-curves:1}.     
   \end{proof}

   \begin{proof}[Proof of Theorem \ref{thm:smooth-add-ruled-surfaces}]
     Let $C\in V_{|D|}^{irr}(\ks_1,\ldots,\ks_r)$. It suffices to show
     that the cohomology group
     $h^1\big(\Sigma,\kj_{X^*(C)/\Sigma}(D)\big)$ vanishes. 

     Suppose this is not the case. Since for $X_0\subseteq X^*(C)$ any
     local complete intersection scheme and
     $z\in\supp(X)$ we have
     \begin{displaymath}
       \deg(X_z)\leq\gamma_0^*(C,z),
     \end{displaymath}
     and since $\gamma\leq \frac{1}{4}$, 
     Lemma \ref{lem:smooth-add-A} applies and there is curve
     $\Delta\sim_a\alpha\cdot C_0+\beta\cdot F$ and a local complete intersection
     scheme $X_0\subseteq X^*(C)$ satisfying the assumptions (a)-(d) there
     and Equation~\eqref{eq:smooth-add-A:1}. 

     Remember that the N\'eron--Severi group of $\Sigma$ is generated
     by a section $C_0$ of $\pi$ and a fibre $F$ with intersection
     pairing given by
     $\left(\begin{smallmatrix}
         -e&1\\1&0 
       \end{smallmatrix}\right)$. 
     Then $K_\Sigma\sim_a-2C_0+(2g-2-e)\cdot F$.
     Note that 
     \begin{displaymath}
       \alpha\geq 0\;\;\;\; \text{ and }\;\;\;\; \beta':=\beta-\tfrac{e}{2}\alpha\geq 0.
     \end{displaymath}
     If we set $b'=b-\frac{ae}{2}$, $\kappa_1=a+2$ and
     $\kappa_2=b+2-2g-\tfrac{ae}{2}=b'+2-2g$, we get
     \begin{equation}\label{eq:smooth-add-ruled-surfaces:3}
       (D-K_\Sigma)^2=-e\cdot(a+2)^2+2\cdot(a+2)\cdot(b+2+e-2g)=2\cdot\kappa_1\cdot\kappa_2.
     \end{equation}

     Fixing the
     notation $\supp(X_0)=\{z_1,\ldots,z_s\}$, $X_i=X_{0,z_i}$,
     and
     $\varepsilon_i=\min\{\deg(X_i),i(C,\Delta;z_i)-\deg(X_i)\}\geq 1$, 
     the conditions on $\Delta$ and $\deg(X_0)$ take the form
     \begin{enumerate}
     \item $a\beta'+b'\alpha\geq \deg(X_0)+\sum_{i=1}^s\varepsilon_i$, 
     \item $\deg(X_0)\geq \kappa_1\cdot \beta'+
       \kappa_2\cdot \alpha-2\alpha\beta'$, and 
     \item $0\leq\alpha\leq \frac{\kappa_1}{2}$ and $0\leq\beta'\leq \frac{\kappa_2}{2}$.
     \end{enumerate}
     The last inequalities follow from (d) in Lemma
     \ref{lem:smooth-add-A} replacing the ample divisor $H$ by the nef
     divisors $F$ respectively $C_0+\frac{e}{2}\cdot F$. 

     From (b) and (c) we deduce
     \begin{displaymath}
       \deg(X_0)\geq \frac{\kappa_1}{2}\cdot \beta'+
       \frac{\kappa_2}{2}\cdot \alpha,
     \end{displaymath}
     and thus, taking \eqref{eq:smooth-add-ruled-surfaces:3}
     into account,
     \begin{equation}
       \label{eq:smooth-add-ruled-surfaces:2}
       \deg(X_0)^2\geq 4\cdot
       \frac{\kappa_1}{2}\cdot\frac{\kappa_2}{2}\cdot \alpha\cdot\beta'
       = \frac{(D-K_\Sigma)^2}{2}\cdot \alpha\cdot \beta'.
     \end{equation}
     Considering now (a) and (b) we get
     \begin{displaymath}
       0<\sum_{i=1}^s\varepsilon_i\leq \Delta.(\Delta+K_\Sigma)=
       2\alpha\beta'+(2g-2)\cdot \alpha -2\beta'
       \leq \frac{\alpha\beta'}{2\gamma},
     \end{displaymath}
     where the last inequality holds if $\beta'\not=0$. We see, in
     particular, that $\beta'\not=0$ if $g\leq 1$.
     But this together with \eqref{eq:smooth-add-ruled-surfaces:2}
     gives for $\beta'\not=0$
     \begin{displaymath}
       \sum_{i=1}^s\varepsilon_i\leq\frac{\deg(X_0)^2}{\gamma\cdot(D-K_\Sigma)^2}.
     \end{displaymath}

     If $\beta'=0$, then we deduce from (a) and (b)
     \begin{displaymath}
       0<\sum_{i=1}^s\varepsilon_i\leq 
       (2g-2)\cdot \alpha\leq
       4\cdot (g-1)\cdot A\cdot \frac{\deg(X_0)^2}{(D-K_\Sigma)^2}
       \leq \frac{\deg(X_0)^2}{\gamma\cdot(D-K_\Sigma)^2}.
     \end{displaymath}

     Applying the Cauchy inequality, we finally get
     \begin{displaymath}
       \gamma\cdot(D-K_\Sigma)^2\leq
       \frac{\deg(X_0)^2}{\sum_{i=1}^s\varepsilon_i}
       \leq \sum_{i=1}^s\frac{\deg(X_i)^2}{\varepsilon_i}
       \leq \sum_{i=1}^r\gamma_0^*(\ks_i),
     \end{displaymath}
     in contradiction to Assumption \eqref{eq:smooth-add-ruled-surfaces:1}.     
   \end{proof}

   \bibliographystyle{amsalpha-tom}
   \bibliography{bibliographie}

\providecommand{\bysame}{\leavevmode\hbox to3em{\hrulefill}\thinspace}
\begin{thebibliography}{Lue87b}

\bibitem[ChC99]{CC99}
Luca Chiantini and Ciro Ciliberto, \emph{On the {Severi} variety of surfaces in
  {$\mathbbm P_{\mathbbm C}^3$}}, J.\ Algebraic Geom. \textbf{8} (1999),
  67--83.

\bibitem[ChS97]{CS97}
Luca Chiantini and Edoardo Sernesi, \emph{Nodal curves on surfaces of general
  type}, Math.\ Ann. \textbf{307} (1997), 41--56.

\bibitem[Fla01]{Fla02}
Flaminio Flamini, \emph{Moduli of nodal curves on smooth surfaces of general
  type}, J.\ Algebraic Geom. \textbf{11} (2001), no.~4, 725--760.

\bibitem[FlM01]{FM01}
Flaminio Flamini and C.\ Madonna, \emph{Geometric linear normality for nodal
  curves on some projective surfaces}, no.~1, 269--283.

\bibitem[GLS97]{GLS97}
{Gert-Martin} Greuel, Christoph Lossen, and Eugenii Shustin, \emph{New
  asymptotics in the geometry of equisingular families of curves}, Internat.\
  Math.\ Res.\ Notices \textbf{13} (1997), 595--611.

\bibitem[GLS98]{GLS98b}
{Gert-Martin} Greuel, Christoph Lossen, and Eugenii Shustin, \emph{Geometry of
  families of nodal curves on the blown-up projective plane}, Trans.\ Amer.\
  Math.\ Soc. \textbf{350} (1998), 251--274.

\bibitem[GLS00]{GLS00}
{Gert-Martin} Greuel, Christoph Lossen, and Eugenii Shustin, \emph{Castelnuovo
  function, zero-dimensional schemes, and singular plane curves}, J.\ Algebraic
  Geom. \textbf{9} (2000), no.~4, 663--710.

\bibitem[GrK89]{GK89}
{Gert-Martin} Greuel and Ulrich Karras, \emph{Families of varieties with
  prescribed singularities}, Comp.\ math. \textbf{69} (1989), 83--110.

\bibitem[GrL96]{GL96}
{Gert-Martin} Greuel and Christoph Lossen, \emph{Equianalytic and equisingular
  families of curves on surfaces}, Manuscripta Math. \textbf{91} (1996),
  323--342.

\bibitem[Har77]{Har77}
Robin Hartshorne, \emph{Algebraic geometry}, Springer, 1977.

\bibitem[Laz97]{Laz97}
Robert Lazarsfeld, \emph{Lectures on linear series}, Complex Algebraic Geometry
  (J\'anos Koll\'ar, ed.), IAS/Park City Mathematics Series, no.~3, Amer.\
  Math.\ Soc., 1997, pp.~161--219.

\bibitem[LoK03]{LK03}
Christoph Lossen and Thomas Keilen, \emph{The $\gamma_\alpha$-invariant for
  plane curve singularities}, Preprint, 2003.

\bibitem[Los98]{Los98}
Christoph Lossen, \emph{The geometry of equisingular and equianalytic families
  of curves on a surface}, Phd thesis, FB Mathematik, Universit\"at
  Kaiserslautern, Aug. 1998, http:// \!\!www. \!\!mathematik. \!\!uni-kl.
  \!\!de/ \!\!\textasciitilde lossen/ \!\!download/ \!\!Lossen002/
  \!\!Lossen002.ps.gz.

\bibitem[Lue87a]{Lue87}
Ignacio Luengo, \emph{The $\mu$-constant stratum is not smooth}, Inventiones
  Math. \textbf{90} (1987), 139--152.

\bibitem[Lue87b]{Lue87a}
Ignacio Luengo, \emph{On the existence of complete families of plane curves,
  which are obstructed}, J.\ LMS \textbf{36} (1987), 33--43.

\bibitem[Sev21]{Sev21}
Francesco Severi, \emph{{Vorlesungen \"uber Algebraische Geometrie}},
  Bibliotheca Mathematica Teubneriana, no.~32, Teubner, 1921.

\bibitem[Shu87]{Shu87}
Eugenii Shustin, \emph{Versal deformation in the space of plane curves of fixed
  degree}, Funct.\ An.\ Appl. \textbf{21} (1987), 82--84.

\bibitem[Shu91]{Shu91a}
Eugenii Shustin, \emph{On manifolds of singular algebraic curves}, Selecta
  Math.\ Sov. \textbf{10} (1991), 27--37.

\bibitem[Shu94]{Shu94}
Eugenii Shustin, \emph{Smoothness and irreducibility of varieties of algebraic
  curves with nodes and cusps}, Bull.\ SMF \textbf{122} (1994), 235--253.

\bibitem[Shu96]{Shu96a}
Eugenii Shustin, \emph{Smoothness and irreducibility of varieties of algebraic
  curves with ordinary singularities}, Israel Math.\ Conf.\ Proc., no.~9,
  Amer.\ Math.\ Soc., 1996, pp.~393--416.

\bibitem[Shu97]{Shu97}
Eugenii Shustin, \emph{Smoothness of equisingular families of plane algebraic
  curves}, Math.\ Res.\ Not. \textbf{2} (1997), 67--82.

\bibitem[Tan80]{Tan80}
Allen Tannenbaum, \emph{Families of algebraic curves with nodes}, Compositio
  Math. \textbf{41} (1980), 107--126.

\bibitem[Tan82]{Tan82}
Allen Tannenbaum, \emph{Families of curves with nodes on {K3}-surfaces}, Math.\
  Ann. \textbf{260} (1982), 239--253.

\bibitem[Vas90]{Vas90}
Victor~A.\ Vassiliev, \emph{Stable cohomology of complements to the
  discriminants of deformations of singularities of smooth functions}, J.\
  Sov.\ Math. \textbf{52} (1990), 3217--3230.

\bibitem[Wah74]{Wah74a}
Jonathan~M.\ Wahl, \emph{Deformations of plane curves with nodes and cusps},
  Amer.\ J.\ Math. \textbf{96} (1974), 529--577.

\end{thebibliography}

\end{document}